\documentclass[12pt,twoside]{amsart}
\usepackage{amssymb}
\usepackage{amscd}

\title[Toric Mori theory]
{Introduction to the toric Mori theory}  
\author{Osamu Fujino and Hiroshi Sato} 
\thanks{The second author is partly supported by the 
Grant-in-Aid for JSPS Fellows, The Ministry of 
Education, Science, Sports and Culture, Japan.}
\address{Graduate School of Mathematics\\ 
 Nagoya University, Chikusa-ku Nagoya 464-8602 Japan}
\email{fujino@math.nagoya-u.ac.jp}
\address{Department of Mathematics\\ 
 Tokyo Institute of Technology, 2-12-1 Oh-Okayama, 
Meguro-ku, Tokyo 152-8551, Japan}
\email{hirosato@math.titech.ac.jp}
\newcommand{\xdiscrep}[0]{{\operatorname{discrep}}}
\newcommand{\xSupp}[0]{{\operatorname{Supp}}}
\newcommand{\Exc}[0]{{\operatorname{Exc}}}
\newcommand{\Pic}[0]{{\operatorname{Pic}}}
\newcommand{\Proj}[0]{{\operatorname{Proj}}}
\newcommand{\Star}[0]{{\operatorname{Star}}}
\newcommand{\Hom}[0]{{\operatorname{Hom}}}
\newcommand{\Sing}[0]{{\operatorname{Sing}}}
\newtheorem{thm}{Theorem}[section]
\newtheorem{lem}[thm]{Lemma}
\newtheorem{cor}[thm]{Corollary}
\newtheorem{prop}[thm]{Proposition}
\newtheorem{lem-def}[thm]{Definition-Lemma}

\theoremstyle{definition}
\newtheorem{ex}[thm]{Example}
\newtheorem{defn}[thm]{Definition}
\newtheorem{rem}[thm]{Remark}
\newtheorem*{ack}{Acknowledgments}       
 
\newtheorem*{notation}{Notation}         

\newtheorem*{claim}{Claim}         

\newtheorem{say}[thm]{}
\begin{document}
\bibliographystyle{amsalpha+}

\begin{abstract}
The main purpose of this paper is to give a 
simple and non-combinatorial proof of the 
toric Mori theory. 
Here, the toric Mori theory means the (log) Minimal 
Model Program (MMP, for short) for toric varieties. 
We minimize the arguments on fans and their decompositions. 
We recommend this paper to the following people: 
\begin{itemize}
\item[(A)] those who are uncomfortable with manipulating fans and 
their decompositions, 
\item[(B)] those who are familiar with toric geometry but not 
with the MMP. 
\end{itemize}
People in the category (A) will be 
relieved from tedious combinatorial arguments in several 
problems. 
Those in the category (B) will discover the potential of the 
toric Mori theory. 

As applications, we treat the Zariski decomposition on 
toric varieties and the partial resolutions of 
non-degenerate hypersurface singularities. 
By these applications, the reader will learn 
to use the toric Mori theory. 
\end{abstract}

\maketitle
\tableofcontents

\section{Introduction}\label{intro}
The main purpose of this paper is to give a 
simple and non-combinatorial proof of the 
toric Mori theory. 
Here, the toric Mori theory means the (log) Minimal 
Model Program (MMP, for short) for toric varieties. 

In his famous and beautiful paper \cite{reid}, Reid 
carried out the toric Mori theory under the assumption that 
the variety is {\em{complete}}. 
His arguments are combinatorial. 
Thus, it is not so obvious whether we can remove the completeness 
assumption from his paper. 
We quote his idea from \cite{reid}. 

\begin{quote}
(0.3) {\em{Remarks.}} The hypothesis that $A$ is complete is not 
essential; it can be reduced to the projective case, or possibly 
eliminated by a careful (and rather tedious) 
rephrasing of the arguments 
of \S\S 1--3. The projectivity hypothesis on 
$f$ is needed in order for the statement of (0.2) to 
make sense, since without projectivity the cone 
${\rm NE}(V/A)$ will usually not have any extremal rays.   
\end{quote}

We prefer not to simply rephrase his approach, 
which entails tedious combinatorial arguments. 
Instead, our proof, which is independent of 
Reid's proof, heavily relies on the general machinery of 
the Minimal Model Program 
and the special properties of toric varieties.  
Thus, our proof works without the completeness assumption. 

For the details of the toric Mori theory, see 
\cite[\S 2.5]{oda}, \cite[\S 5.2]
{kmm}, \cite{op}, \cite{lat}, \cite{shihoko2}, 
\cite[Chapter 14]{ma}, \cite{wis}, and 
\cite{fujino}. 
Matsuki \cite{ma} corrected some minor errors in 
\cite{reid} and pointed out some ambiguities in \cite{reid} 
and \cite{kmm}. 
See Remarks 14-1-3 (ii), 14-2-3, and 14-2-7 in \cite{ma}.  
We believe that these remarks help the reader to understand 
\cite{reid}. 
We recommend that the reader compare this paper with 
\cite[Chapter 14]{ma}. 
Shokurov \cite{shokurov} treats the MMP 
for toric varieties in a non-combinatorial way (see 
\cite[Example 3]{shokurov}). 
His arguments are quite different from ours. 
For the more advanced topics of the toric 
Mori theory, see \cite{fujino2}. 

For the outline of the general MMP, which is 
still conjectural in dimension $\geq 4$, 
see \cite[Introduction]{kmm} 
or \cite[2.14, \S3.7]{km}. 

We note that the Zariski decomposition on 
toric varieties has already been treated by various researchers. 
The reason we treat it here is to show that the 
Zariski decomposition on toric varieties is an 
easy consequence of the toric Mori theory. 
The partial resolutions of 
non-degenerate hypersurface singularities were 
treated by Ishii, who 
divided the cone by the data 
of the Newton polytope. 
One of her proofs contains a gap 
(see Remark \ref{shihotan}), 
but the results themselves are correct. 
Her results also become easy consequences of 
the toric Mori theory. 

Note that we cannot recover combinatorial 
aspects of \cite{reid} and \cite{shihoko} by 
our method. 
So, this paper does not depreciate \cite{reid} and 
\cite{shihoko}. 

Almost all the results in this paper are more or less known 
to the experts. However, 
some of them were not stated explicitly before. 
Some of the proofs that we give in this paper are 
new and much simpler than the known ones. 
We hope that this  paper will help the reader to understand 
the toric Mori theory. 


\vspace{5mm}

We summarize the contents of this paper: 
In Section \ref{pre}, we fix the notation and collect basic 
results. 
Section \ref{sec1} explains the toric Mori theory. 
Section \ref{proofof} is the main part of this paper. 
There, we give a simple and non-combinatorial 
proof of the toric Mori theory. 
In Section \ref{sec-zari}, we consider the Zariski decomposition 
on toric varieties. 
In Section \ref{appl}, we apply the toric Mori theory 
to the study of the partial resolutions of 
non-degenerate hypersurface singularities. 
We reprove Ishii's results. 

\begin{notation} 
Here is a list of some of the standard notation we use. 
\begin{itemize}
\item[(1)] For a real number $d$, its {\em{round down}} is the 
largest integer $\leq d$. 
It is denoted by $\llcorner d\lrcorner$. 
If $D=\sum d_iD_i$ is a divisor with real coefficients and  the 
$D_i$ are distinct prime divisors, then 
we define the {\em{round down}} of $D$ as 
$\llcorner D\lrcorner:=\sum \llcorner d_i\lrcorner D_i$. 
\item[(2)] Let $f:X\longrightarrow Y$ be a proper birational 
morphism between normal varieties. 
Then, $f$ is said to be {\em{small}} 
if $f$ is an isomorphism in codimension one.  
\item[(3)] The symbol 
$\mathbb Z_{\geq 0}$ (resp.~$\mathbb Q_{\geq 0}$, 
$\mathbb R_{\geq 0}$) 
denotes the set of non-negative integers 
(resp.~rational numbers, real numbers). 
\end{itemize}
We will work over an algebraically closed field $k$ 
throughout this paper. 
The characteristic of $k$ is arbitrary from Section \ref{pre} to 
Section \ref{sec-zari} unless otherwise stated. 
In Section \ref{appl}, we assume that $k$ is the complex 
number field $\mathbb C$. 
\end{notation}

\begin{ack} 
We would like to thank Professors 
Shigefumi Mori and Tatsuhiro Minagawa for fruitful discussions 
and useful comments. 
We are grateful to Professor Shihoko Ishii, who 
kindly answered our questions about the proof of 
\cite[Theorem 3.1]{shihoko} and 
gave us several useful comments. 
Finally, we thank Professor Tadao Oda, whose 
comments on the preliminary version made this 
paper more readable. 
\end{ack}

\section{Preliminaries}\label{pre}

\subsection{Toric varieties} 
In this subsection, we recall the basic notion of toric varieties 
and fix the notation. 
For the basic results about toric varieties, see 
\cite{kmms}, \cite{danilov}, \cite{oda}, and \cite{fulton}. 

\begin{say} 
Let $N\simeq \mathbb Z^n$ be a lattice of rank $n$. 
A toric variety $X(\Delta)$ is associated to a {\em{fan}} $\Delta$, 
a finite collection of convex cones $\sigma\subset N_\mathbb R :=
N\otimes _{\mathbb Z}\mathbb R$ satisfying the following: 
\begin{enumerate}
\item[(i)] 
Each convex cone $\sigma\in\Delta$ is rational 
polyhedral in the sense 
that there are finitely many 
$v_1, \cdots, v_s\in N\subset N_{\mathbb R}$ such 
that 
$$
\sigma=\{r_1v_1+\cdots +r_sv_s; \ r_i\in\mathbb R_{\geq 0}
\mbox{ for all }i\}
$$ 
and is strongly convex in the sense that
$$
\sigma \cap (-\sigma)=\{0\}. 
$$
\item[(ii)] Each face $\tau$ of a convex cone $\sigma\in \Delta$ 
again belongs to $\Delta$. 
\item[(iii)] The intersection of two cones in $\Delta$ is a face of 
each. 
\end{enumerate}

\begin{defn}\label{saisyo}
The {\em{dimension}} $\dim \sigma$ of $\sigma$ is 
the dimension of the linear space 
$\mathbb R\cdot \sigma=\sigma +(-\sigma)$ spanned 
by $\sigma$. 
\end{defn}

We define the sublattice $N_{\sigma}$ of 
$N$ generated (as a subgroup) by $\sigma\cap N$ as 
follows: 
$$
N_{\sigma}:=\sigma\cap N+(-\sigma\cap N). 
$$
\end{say}

The {\em{star}} of a cone $\tau$ can be defined abstractly 
as the set of cones $\sigma$ in $\Delta$ that 
contain $\tau$ as a face. For any such cone $\sigma$, the image 
$$
\overline \sigma=\left( \sigma+(N_{\tau})_{\mathbb R}
\right) / 
(N_{\tau})_{\mathbb R}\subset N(\tau)_{\mathbb R}
$$ 
by the projection $N\to N(\tau):=N/{N_{\tau}}$ 
is a cone in $N(\tau)$. 
These cones $\{\overline \sigma ; \tau\prec \sigma\}$ 
form a fan in $N(\tau)$, and we denote this fan by 
$\Star(\tau)$. 
We set $V(\tau)=X(\Star (\tau))$. 
It is well-known that $V(\tau)$ is an $(n-k)$-dimensional 
closed toric subvariety of $X(\Delta)$, where $\dim \tau=k$. 
If $\dim V(\tau)=1$ (resp.~$n-1$), then we call $V(\tau)$ 
a {\em{torus invariant curve}} (resp.~{\em{torus invariant 
divisor}}). 
For details on the correspondence between $\tau$ and 
$V(\tau)$, see, for instance, \cite[3.1 Orbits]{fulton}. 

\begin{defn}[{$\mathbb Q$-Cartier divisor and 
$\mathbb Q$-factoriality}]\label{q-fac}
Let $D=\sum d_iD_i$ be a $\mathbb Q$-divisor on 
a normal variety $X$, that is, 
$d_i\in \mathbb Q$ and $D_i$ is a prime divisor 
on $X$ for every $i$. 
Then $D$ is {\em{$\mathbb Q$-Cartier}} if 
there exists a positive integer $m$ such that 
$mD$ is a Cartier divisor. 
A normal variety $X$ is said to be {\em{$\mathbb Q$-factorial}} 
if every prime divisor $D$ on $X$ is $\mathbb Q$-Cartier. 
\end{defn}

The next lemma is well-known. 
See, for example, \cite[(1.9)]{reid} or 
\cite[Lemma 14-1-1]{ma}. 

\begin{lem}\label{kantan}
A toric variety $X(\Delta)$ is $\mathbb Q$-factorial 
if and only if each cone $\sigma\in \Delta$ is simplicial. 
\end{lem}

The following remarks are easy but important. 

\begin{rem}
Let $D$ be a Cartier (resp.~$\mathbb Q$-Cartier) 
divisor on a toric variety. 
Then $D$ is linearly (resp.~$\mathbb Q$-linearly) 
equivalent to 
a torus invariant divisor (resp.~$\mathbb Q$-divisor). 
\end{rem}

\begin{rem}[{cf.~\cite[(4.1)]{reid}}]
Let $X$ be a toric variety and $D$ the complement of the 
big torus regarded as a reduced divisor. Then $K_X+D\sim 0$. 
\end{rem}

\begin{say}[the Kleiman-Mori Cone]
Let $f:X\longrightarrow Y$ be a proper morphism between 
normal varieties $X$ and $Y$; 
a $1$-cycle of $X/Y$ is a formal sum $\sum a_iC_i$ with 
complete curves $C_i$ in the fibers of $f$, and 
$a_i\in \mathbb Z$. 
We put 
$$
Z_1(X/Y):=\{1\text{-cycles of} \ X/Y\}, 
$$
and 
$$
Z_1(X/Y)_{\mathbb Q}:=
Z_1(X/Y)\otimes \mathbb Q.  
$$
There is a pairing 
$$
\Pic (X)\times Z_1(X/Y)_{\mathbb Q}
\longrightarrow \mathbb Q
$$
defined by 
$(\mathcal L, C)\mapsto \deg _C\mathcal L$, extended 
by bilinearity. 
Define 
$$
N^1(X/Y):=(\Pic (X)\otimes \mathbb Q)/\equiv
$$
and 
$$
N_1(X/Y):= Z_1(X/Y)_{\mathbb Q}/\equiv, 
$$ 
where the {\em numerical equivalence} $\equiv$ 
is by definition the smallest equivalence relation which 
makes $N^1$ and $N_1$ into dual 
spaces. 

Inside $N_1(X/Y)$ there is a distinguished cone of effective 
$1$-cycles, 
$$
{\rm NE}(X/Y)=\{\, Z\, | \ Z\equiv 
\sum a_iC_i \ \text{with}\ a_i\in \mathbb Q_{\geq 0}\}
\subset N_1(X/Y).  
$$
A subcone $F\subset {\rm NE}(X/Y)$ is said to be {\em{extremal}} 
if $u,v\in {\rm NE}(X/Y)$, $u+v\in F$ imply $u,v\in F$. 
The cone $F$ is also called an 
{\em{extremal face}} of ${\rm NE}(X/Y)$. 
A one-dimensional extremal face is called an {\em{extremal 
ray}}. 

We define the {\em{relative Picard number}} $\rho(X/Y)$ by
$$
\rho (X/Y):=\dim _{\mathbb Q}N^1(X/Y)< \infty. 
$$
An element $D\in N^1(X/Y)$ is called {\em{$f$-nef}} if 
$D\geq 0$ on ${\rm NE}(X/Y)$. 

If $X$ is complete and $Y$ is a point, 
write ${\rm NE}(X)$ and $\rho(X)$ 
for ${\rm NE}(X/Y)$ and $\rho(X/Y)$, respectively. 
We note that $N_1(X/Y)\subset N_1(X)$, and $N^1(X/Y)$ 
is the corresponding 
quotient of $N^1(X)$. 
\end{say}

\subsection{Singularities of pairs}

In this subsection, we quickly review the definitions of singularities 
which we use in the MMP. 
For details, see, for example, \cite[\S 2.3]{km}. 
We recommend that the reader skip this subsection on 
first reading. 

\begin{say} 
Let us recall the definitions of the singularities for pairs. 

\begin{defn}[Discrepancies and Singularities of Pairs]\label{def1}
Let $X$ be a normal variety and $D=\sum d_i D_i$ a 
$\mathbb Q$-divisor on $X$, 
where $D_i$ are distinct and irreducible 
such that $K_X+D$ is $\mathbb 
Q$-Cartier. 
Let $f:Y\longrightarrow X$ be a proper birational morphism 
from a normal variety $Y$. 
Then we can write 
$$
K_Y=f^{*}(K_X+D)+\sum a(E,X,D)E, 
$$ 
where the sum runs over all the distinct prime divisors $E\subset Y$, 
and $a(E,X,D)\in \mathbb Q$. This $a(E,X,D)$ is called the 
{\em discrepancy} of $E$ with respect to $(X,D)$. 
We define 
$$
\xdiscrep (X,D):=\inf _{E}\{a(E,X,D)\, 
|\, E \text{ is exceptional over}\  X  \}.
$$ 
From now on, we assume 
that $0\leq d_i \leq 1$ 
for every $i$. 
We say that $(X,D)$ is 
$$
\begin{cases}
\text{terminal}\\
\text{canonical}\\
\text{klt}\\
\text{plt}\\
\text{lc}\\
\end{cases}
\quad {\text{if}} \quad \xdiscrep (X,D) 
 \quad
\begin{cases}
>0,\\
\geq 0,\\
>-1\quad {\text {and \quad $\llcorner D\lrcorner =0$,}}\\
>-1,\\
\geq -1.\\
\end{cases}
$$ 
Here klt is an abbreviation for 
{\em {Kawamata log terminal}}, plt 
for {\em{purely log terminal}}, and lc for 
{\em{log-canonical}}. 

If there exists a log resolution $f:Y\longrightarrow X$ 
of $(X,D)$, that is, $Y$ is non-singular, the exceptional 
locus $\Exc(f)$ is a divisor, and $\Exc(f)\cup 
f^{-1}(\xSupp D)$ is a simple normal crossing 
divisor, such 
that $a(E_i,X,D)>-1$ for every exceptional divisor $E_i$ on 
$Y$, then the pair $(X,D)$ is said to be {\rm{dlt}}, 
which is an abbreviation for {\em{divisorial log terminal}}. 
\end{defn}

For details on dlt, see \cite[Definition 2.37, 
Theorem 2.44]{km}. 
The following results are well-known to the experts. 
See, for example, \cite[Lemma 5.2]{fujino} or 
\cite[Proposition 14-3-2]{ma}.  

\begin{prop}\label{ata} 
Let $X$ be a toric variety and $D$ the complement of the 
big torus regarded as a reduced divisor. 
Then $(X,D)$ is log-canonical. 
Let $D=\sum _iD_i$ be the irreducible decomposition of $D$. 
We assume that $K_X+\sum _i a_iD_i$ is $\mathbb Q$-Cartier, 
where $0\leq a_i<1$ $($resp.~$0\leq a_i\leq 1$$)$ 
for every $i$. 
Then $(X,\sum a_iD_i)$ is {\rm klt} 
$($resp.~log-canonical$)$.  
\end{prop}
\end{say}

\section{Toric Mori theory}\label{sec1}

We will work over an algebraically closed field $k$ 
of arbitrary characteristic throughout this section. 
Let us explain the Minimal Model Program for toric varieties. 

\begin{say}[Minimal Model Program for Toric Varieties]\label{mmp}
We start with a projective 
toric morphism $f:X\longrightarrow Y$, 
that is, $f$ is induced by a map of lattices, 
where $X=:X_0$ is a $\mathbb Q$-factorial toric variety, 
and a $\mathbb Q$-divisor 
$D_0:=D$ on $X$. 
The aim is to set up a recursive procedure that 
creates intermediate $f_i:X_i\longrightarrow Y$ and $D_i$ 
on $X_i$. 
After finitely many steps, we obtain final objects 
$\widetilde {f}:\widetilde {X}\longrightarrow Y$ and 
$\widetilde {D}$ on $\widetilde{X}$. 
Assume that we have already constructed 
$f_i:X_i\longrightarrow Y$ 
and $D_i$ with 
the following properties: 
\begin{enumerate}
\item[(i)] $X_i$ is $\mathbb Q$-factorial and 
$f_i$ is projective, 
\item[(ii)] $D_i$ is a $\mathbb Q$-divisor on $X_i$. 
\end{enumerate}

If $D_i$ is $f_i$-nef, then we set $\widetilde{X}:=X_i$ 
and $\widetilde {D}=D_i$. 
If $D_i$ is not $f_i$-nef, 
then we can take an 
extremal ray $R$ of ${\rm NE}(X_i/Y)$ 
such that $R\cdot D_i<0$ (see 
Theorem \ref{cone} below). 
Thus we have a contraction morphism $\varphi_R:X_i
\longrightarrow W_i$ over 
$Y$ (see Theorem \ref{cont} below). 
If $\dim W_i <\dim X_i$ (in which case we call $\varphi_R$ a 
{\em{Fano contraction}}), then we set $\widetilde{X} 
:=X_i$ and $\widetilde{D}:=D_i$ 
and stop the process. 
If $\varphi_R$ is birational and contracts a divisor (we call this 
a {\em{divisorial contraction}}), 
then we put $X_{i+1}:=W_i$, $D_{i+1}:=\varphi_{R*}D_i$ and 
repeat this process. 
In the case where $\varphi_R$ is small 
(we call this a {\em{flipping 
contraction}}), then there 
exists a log-flip $\psi:X_i\dasharrow X_i^{+}$ over $Y$. 
Here, a log-flip means an {\em{elementary transformation}} 
with respect to $D_i$ (see Theorem \ref{flip} below). 
Note that $\psi$ is an isomorphism in codimension one. 
We put $X_{i+1}:=X_i^{+}$, 
$D_{i+1}:=\psi_{*}D_i$ and repeat this process. 
By counting the relative Picard number $\rho(X_i/Y)$, 
divisorial contractions can occur finitely many times 
(see Theorem \ref{cont}). 
By Theorem \ref{termi} below, every sequence of 
log-flips terminates after finitely many steps. 
So, this process always terminates and 
we obtain $\widetilde{f}:\widetilde{X}\longrightarrow Y$ 
and $\widetilde{D}$. 
We call this process the ($D$-){\em{Minimal Model Program 
over}} $Y$, 
where $D$ is a divisor used in the process. 
When we apply the Minimal Model Program (MMP, for short), 
we say that, 
for example, we {\em{run the MMP over $Y$ with respect to the 
divisor $D$}}. 
\end{say}

\begin{rem}[Toric Mori theory vs the general MMP] 
The general (log) MMP is still conjectural in dimension 
$\geq 4$ (see \cite[\S 3.7]{km}). 
As we see in Section \ref{proofof}, the MMP for 
toric varieties is fully established and it works 
for {\em{any}} $\mathbb Q$-divisor $D$. 
A generalization of the MMP 
for non-$\mathbb Q$-factorial toric varieties is treated in 
\cite[Section 2]{fujino2}. 
\end{rem}

\begin{rem}[Original toric Mori theory]  
In the above MMP, if we assume that $Y$ is projective, 
$X$ has only terminal singularities, $f$ is birational, and 
$D=K_X$ the canonical divisor of $X$, 
then we can recover the original toric Mori theory 
in \cite[(0.2) Theorem]{reid}.  
\end{rem}

\begin{rem}
In \ref{mmp}, it is sufficient to assume that 
$D$ is an $\mathbb R$-divisor. 
We do not treat the $\mathbb R$-generalization here 
because this generalization is 
obvious for experts and we do not need $\mathbb R$-divisors 
in this paper. 
We leave the details to the reader. 
\end{rem}

\begin{rem}\label{qfac}
In general, the assumption that 
$X$ is $\mathbb Q$-factorial is not 
crucial in the toric Mori theory. 
By \cite[Corollary 5.9]{fujino}, 
there always exists a small projective 
$\mathbb Q$-factorialization. Namely, for any toric 
variety $X$, there exists a small projective toric 
morphism $\widetilde{X}\longrightarrow X$ 
such that $\widetilde{X}$ 
is $\mathbb Q$-factorial. 
We refer the reader to \cite[Section 5]{fujino}, 
which is a baby version of this paper. 
We note that 
we can remove the assumption that $X$ is complete in \cite
[Theorem 5.5, Proposition 5.7]{fujino} by the results in 
Sections \ref{sec1} and \ref{proofof} of this paper. 
\end{rem}

\begin{rem}\label{qcar}
Let $X$ be a toric variety and $D$ a $\mathbb Q$-divisor on $X$. 
The assumption that $D$ is $\mathbb Q$-Cartier 
can be removed in some cases. 
In fact, by replacing $X$ by its small projective 
$\mathbb Q$-factorialization we can 
assume that $D$ is $\mathbb Q$-Cartier. 
See the proof of Corollary \ref{finifini}. 
\end{rem}

\section{Proof of the toric Mori theory}\label{proofof}

In this section, we give a simple and non-combinatorial 
proof of the toric Mori theory. 

\begin{thm}[the Cone Theorem]\label{cone}
Let $f:X\longrightarrow Y$ be a 
proper toric morphism. 
Then the cone 
$$
{\rm NE}(X/Y)\subset N_1(X/Y) 
$$ 
is a polyhedral convex cone. 
Moreover, if $f$ is projective, then the cone is strongly convex. 
\end{thm}

\begin{proof} 
By taking the Stein factorization of $f$, we may assume that $f$ is 
surjective with connected fibers. 
We consider $V(\sigma)\subset Y$ for 
some cone $\sigma$. 
Then $f^{-1}(V(\sigma))$ is a union of $V(\tau)\subset X$ for 
some cones $\tau$ since $f$ is a proper toric morphism. 
We divide $Y$ into a finite disjoint union of tori $Y=
\amalg_i Y_i$. 
We put $V_i:=f^{-1}(Y_i)\longrightarrow Y_i$. 
Let $\amalg_jV_{ij}$ be the normalization of $V_i$. 
Then we can check 
that $V_{ij}$ is a toric variety for every $i,j$ by 
using the above fact, that is, 
$f^{-1}(V(\sigma))$ is a union of orbit closures, 
inductively on $\dim V(\sigma)$. 
We note that $V_{ij}$ is dominant onto $Y_i$ for every $i,j$ since 
$Y_i$ is a torus. 
So, we obtain a collection of proper surjective toric morphisms 
with connected fibers:~$\{V_{ij}\longrightarrow Y_i\}_{i,j}$. 
By changing the notation $V_{ij}$, 
we write $\{f_i:X_i\longrightarrow Y_i\}_{i}$ for 
$\{V_{ij}\longrightarrow Y_i\}_{i,j}$. 
We note that $i\ne i'$ does not imply $Y_i\ne Y_{i'}$ 
in this notation. 
Since $Y_i$ is a torus, 
$X_i\simeq F_i\times Y_i$ for every $i$, where 
$F_i$ is a complete toric variety 
(cf.~\cite[p.41, Excercise]{fulton}).  

\begin{claim} 
We have the following commutative diagram: 
$$
\begin{matrix}
N_1(F_i)&\simeq &N_1(X_i/Y_i)\\
\cup& &\cup\\
{\rm NE}(F_i)&\simeq &{\rm NE}(X_i/Y_i)
\end{matrix}
$$ 
for every $i$. 
In particular, ${\rm NE}(X_i/Y_i)$ is a polyhedral convex cone for 
every $i$. 
\end{claim}
\begin{proof}[Proof of Claim] 
We consider the cycle map 
$Z_1(F_i)\longrightarrow Z_1(X_i/Y_i)$ 
which is induced by the inclusion $F_i\simeq 
F_i\times\{{\text{a point of}} \ Y_i\}\subset F_i\times 
Y_i\simeq X_i$. 
It induces 
$$
\varphi_i:N_1(F_i)\longrightarrow N_1(X_i/Y_i). 
$$ 
Let $0\ne v \in N_1(F_i)$. Then 
there exists $\mathcal L\in \Pic (F_i)$ such 
that $\mathcal L\cdot  v\ne 0$. 
Let $p_i:X_i\longrightarrow F_i$ be the 
first projection. 
Then ${p_i}^{*}\mathcal L\cdot 
\varphi_i(v)=\mathcal L\cdot v
\ne 0$ by the 
projection formula. 
Therefore, $\varphi _i$ is injective. 
Since $Y_i$ is a torus and $X_i\simeq F_i\times Y_i$, 
it is obvious that $\varphi_i$ is surjective. 
Since ${\rm NE}(F_i)$ is well-known to be 
a polyhedral convex cone 
(cf.~\cite[p.96 Proposition]{fulton}), 
the other parts are obvious. 
\end{proof}
We consider the following commutative diagram: 
$$
\begin{matrix}
\bigoplus _iN_1(X_i/Y_i)&\twoheadrightarrow &N_1(X/Y)\\
\cup&&\cup\\
\bigoplus_i{\rm NE}(X_i/Y_i)&\twoheadrightarrow &{\rm NE}(X/Y). 
\end{matrix}
$$ 
We note that $\bigoplus _iZ_1(X_i/Y_i)\longrightarrow Z_1(X/Y)$ 
is surjective. 
So, by combining it with the previous claim, 
we obtain the 
required cone theorem for ${\rm NE}(X/Y)\subset N_1(X/Y)$. 
The last part follows from Kleiman's criterion.
\end{proof}

We give another proof of Theorem \ref{cone}, 
which works under the 
assumption that the characteristic of $k$ is zero. 
This assumption is required only in 
the following proof. 

\begin{proof}[Proof of {\em{Theorem \ref{cone}}} 
in characteristic zero]
Assume that the characteristic of $k$ is zero. 
First, we further assume that 
$X$ is quasi-projective and $\mathbb Q$-factorial. 
Let $T$ be the big torus of $X$. We put 
$D=\sum _iD_i=X\setminus T$ 
regarded as a reduced divisor. 
We can take an ample $\mathbb Q$-divisor 
$L=\sum_i a_iD_i$ with $0<a_i<1$ for 
every $i$. 
Then 
$$
-(K_X+\sum _i(1-a_i)D_i)\sim \sum _ia_iD_i
$$ 
is ample (obviously, $f$-ample) 
and $(X, \sum _i(1-a_i)D_i)$ is klt by Proposition \ref{ata}. 
So, the well-known relative cone theorem (see, for 
example, \cite[Theorem 4-2-1]{kmm} or 
\cite[Theorem 3.25]{km}) implies that 
${\rm NE}(X/Y)$ is a rational polyhedral convex cone. 
Next, by Chow's lemma and the desingularization theorem, 
we can take a proper birational toric morphism $X'\longrightarrow 
X$ from a non-singular quasi-projective toric variety $X'$. 
So, the general case follows from the above special case. 
Details are left to the reader. 
\end{proof}

\begin{rem}
In Theorem \ref{cone}, if $X$ is complete, then every 
extremal ray of ${\rm NE}(X/Y)$ is spanned by 
torus invariant curves on $X$. Related topics 
are treated in the 
first author's paper \cite{fujino}. 
If $X$ is not complete, then 
${\rm NE}(X/Y)$ is not necessarily spanned by 
a torus invariant curve, 
as the following example shows:
\end{rem}

\begin{ex}
Let $Y$ be a one-dimensional 
(not necessarily complete) 
toric variety. We put $X=Y\times \mathbb P^1$. 
Let $f:X\longrightarrow Y$ be the first projection. 
Then ${\rm NE}(X/Y)$ is a half line. 
When $Y$ is a one-dimensional torus, 
there are no torus invariant 
curves in the fibers of $f$. 
If $Y\simeq \mathbb P^1$ or $\mathbb A^1$, 
then ${\rm NE}(X/Y)$ is spanned by a torus invariant 
curve in a fiber of $f$. 
\end{ex}

The following remark is obvious since 
a torus is a connected linear algebraic group 
(cf.~\cite[Lemma 5]{sumihiro}). 
We present it for the reader's convenience. 

\begin{rem}
Let $T$ be the big torus of $X$. 
Then $T$ acts on ${\rm NE}(X/Y)$. Let $R$ be an 
extremal ray of ${\rm NE}(X/Y)$. 
Then there exists 
a nef torus invariant Cartier divisor $D$ on $X$ 
such that $D\cdot [C]=0$ if and only if $[C]\in R$. 
So, $R$ is $T$-invariant. Thus, 
$T$ acts on ${\rm NE}(X/Y)$ trivially. 
Hence, the action of $T$ on $N_1(X/Y)$ is trivial as well. 
\end{rem}

Therefore, an extremal ray $R$ of ${\rm NE}(X/Y)$ does not 
necessarily contain a torus invariant curve but 
is torus invariant. 

\begin{thm}[the Contraction Theorem]\label{cont}
Let $f:X\longrightarrow Y$ be a projective toric morphism. 
Let $F$ be an extremal face of ${\rm NE}(X/Y)$. 
Then there exists a projective surjective toric morphism 
$$
\varphi_F:X\longrightarrow Z
$$ 
over $Y$ 
with the following properties{\em:} 
\begin{enumerate}
\item[(i)] $Z$ is a toric variety that is 
projective over $Y$. 
\item[(ii)] ${\varphi_{F}}$ has connected fibers. 
\item[(iii)] Let $C$ be a curve in a fiber of $f$. Then $[C]\in F$ 
if and only if $\varphi_{F}(C)$ is a point. 
\end{enumerate}
Furthermore, if $F$ is an extremal ray $R$ and 
$X$ is $\mathbb Q$-factorial, 
then $Z$ is $\mathbb Q$-factorial and $\rho(Z/Y)=\rho(X/Y)-1$ 
if $\varphi_R$ is not small. 
\end{thm}
\begin{proof}
Since ${\rm NE}(X/Y)$ is a 
polyhedral convex cone, we can take an $f$-nef 
Cartier divisor $D$ 
such that $D\cdot [C]\geq 0$ for every $[C]\in {\rm NE}(X/Y)$ 
and $D\cdot [C]=0$ if and only if $[C]\in F$. 
By replacing $D$ by a linearly equivalent divisor, 
we may assume that 
$D$ is a torus invariant Cartier divisor on $X$. 
We put $\varphi_F:X\longrightarrow 
Z$ as in the 
the following 
Proposition \ref{rel}. 
Then $\varphi_F:X\longrightarrow Z$ has the 
required properties. 
The latter part is well-known. See, for example, 
\cite[Proposition 3.36]{km}. 
\end{proof}

\begin{prop}\label{rel}
Let $f:X\longrightarrow Y$ be a proper surjective 
toric morphism 
between toric varieties. 
Let $D$ be an $f$-nef torus invariant Cartier divisor on $X$. 
Then $D$ is $f$-free, that is, $f^*f_*\mathcal O_X(D)
\longrightarrow \mathcal O_X(D)$ is surjective.  
Moreover, we have a projective toric morphism 
$\varphi:X\longrightarrow Z$ 
over $Y$ such that 
\begin{enumerate} 
\item[(i)] $\varphi$ has only connected fibers, and 
\item[(ii)] For any irreducible curve on $X$ with $f(C)$ being 
a point, $\varphi(C)$ is a point if and only if $D\cdot C=0$. 
\end{enumerate}
\end{prop}
\begin{proof}
Let 
$$
\begin{CD}
f:X@>{g}>>\widetilde Y@>{h}>>Y 
\end{CD} 
$$ 
be the Stein factorization of $f$. 
For the first part, we may assume that 
$Y$ is affine (hence so is $\widetilde Y$). 
It is sufficient to prove that $D$ is $g$-free. 
We can apply the argument in 
\cite[p.68, Proposition]{fulton} with minor modifications. 
See also \cite[Chapter IV.~1.8.~Lemma (2)]{nakayama} 
and \cite[Lemma 14-1-11]{ma}. 
For the second part, $\varphi:X\longrightarrow Z:=
\Proj_{\widetilde Y}\bigoplus_{m\geq 0}
g_*\mathcal O_X(mD)$ is equivariant by construction. 
When $\widetilde Y$ 
is a point, it is well-known that $Z$ is a projective 
toric variety 
that is constructed from a suitable polytope. 
Let $T\subset \widetilde Y$ be the big torus. 
Then $g^{-1}(T)\simeq T\times 
F$ for some complete toric variety $F$. 
So, $Z$ contains a torus as 
a non-empty Zariski open set by 
the above case where $\widetilde Y$ is a point. 
It is obvious that $Z$ is normal and has a suitable 
torus action 
by construction. Therefore, $Z$ is the required toric variety.
\end{proof}

\begin{thm}[Finitely Generatedness of Divisorial Algebra]\label{fini} 
Let $f:X\longrightarrow Y$ be a proper 
birational toric morphism and 
$D$ a torus invariant Cartier divisor on $X$. 
Then 
$$
\bigoplus _{m\geq 0}f_*\mathcal O_X(mD)
$$ 
is 
a finitely generated $\mathcal O_Y$-algebra. 
\end{thm}
\begin{proof} 
We may assume that $Y$ is affine. 
So, it is sufficient to show that 
$$\bigoplus _{m\geq 0}H^0(X,\mathcal O_X(mD))$$ 
is a finitely generated $k$-algebra. 
We put $X=X(\Delta)$, that is, $X$ is a toric variety associated 
to a fan $\Delta$ in $N$. 
Let $\psi_D$ be the support function of $D$. 
We put 
$$
P_D=\{u\in M_{\mathbb R} \ |\  u\geq \psi _D \ \text{on} \ 
|\Delta|\}, 
$$
and 
$$
P_{aD}=\{u\in M_{\mathbb R} \ |\ u\geq a\psi _D \ \text{on} \ 
|\Delta|\},  
$$ 
where $M:=\Hom_{\mathbb Z}(N,\mathbb Z)$ is the dual lattice 
of $N$ and $|\Delta|$ stands for the 
{\em{support}} of the fan $\Delta$. 
We define 
$$
C=\{(u,a)\in M_{\mathbb R}\times {\mathbb R}_{\geq 0} \ | \ 
u \in P_{aD}\}, 
$$ 
and $C_{\mathbb Z}=C\cap (M\times {\mathbb Z}_{\geq 0})$. 
The $k$-algebra 
$\bigoplus _{m\geq 0}H^0(X,\mathcal O_X(mD))$ is the 
semi-group ring associated to $C_{\mathbb Z}$. 
We can easily check that $C$ is 
a finite intersection of half spaces, which are 
defined over the rational numbers, 
in $M_{\mathbb R}\times {\mathbb R}$. 
Therefore, $C$ is a rational polyhedral convex cone. 
Thus, $C_{\mathbb Z}$ is a finitely generated semi-group. 
This implies that 
$\bigoplus _{m\geq 0}H^0(X,\mathcal O_X(mD))$ is 
a finitely generated $k$-algebra.  
\end{proof}

We will generalize the above 
theorem in Corollary \ref{finifini} below. 

\begin{thm}[Elementary Transformation]\label{flip}
Let $\varphi:X\longrightarrow W$ be a small toric morphism 
and $D$ a 
torus invariant $\mathbb Q$-Cartier 
divisor on $X$ such that $-D$ is 
$\varphi$-ample. Let $l$ be a positive integer such that 
$lD$ is Cartier. 
Then there exists a small projective toric morphism 
$$
\varphi ^{+}:X^{+}:=
\Proj_W\bigoplus_
{m \geq 0}\varphi_*\mathcal O_X(mlD)\longrightarrow W
$$ 
such that $D^+$ is a $\varphi^+$-ample 
$\mathbb Q$-Cartier divisor, 
where $D^+$ is the proper transform of $D$ on $X^+$. 
The commutative diagram
$$
\begin{matrix}
X & \dashrightarrow & \ X^+ \\
{\ \ \ \ \searrow} & \ &  {\swarrow}\ \ \ \ \\
 \ & W &  
\end{matrix}
$$ 
is called the {\em{elementary transformation (with respect to 
$D$)}}. 

Moreover, if $X$ is $\mathbb Q$-factorial and 
$\rho(X/W)=1$, then so is $X^+$ and $\rho (X^+/W)=1$. 
\end{thm}
\begin{proof}The 
first part is obvious by the above theorem and the construction 
of $\varphi^+:X\longrightarrow W$. 
See, for example, \cite[4.2 Proposition]{FA} or 
\cite[Lemma 6.2]{km}. 
The latter part is well-known. 
See, for example, \cite[Proposition 3.37]{km}. 
\end{proof}

\begin{thm}[Termination of Elementary Transformations]\label{termi} 
Let 
$$
\begin{matrix} 
X_0 & \dashrightarrow & X_1 & 
\dashrightarrow & X_2 
&\dashrightarrow\cdots\\
{\ \ \ \ \searrow} & \ &  {\swarrow}\ \  {\searrow} & \ &  
{\swarrow\ \ \ } &\\
 \ & W_0 & \  \ & W_1 & \ & & 
\end{matrix}
$$
be a sequence of elementary transformations of toric varieties 
with respect to a fixed $\mathbb Q$-Cartier divisor $D$. 
More precisely, the commutative diagram
$$
\begin{matrix}
X_i & \dashrightarrow & \ X_{i+1} \\
{\ \ \ \ \searrow} & \ &  {\swarrow}\ \ \ \ \\
 \ & W_i &  
\end{matrix}
$$ 
is the elementary transformation with respect to 
$D_i$, 
where $D_0=D$ and $D_i$ is the proper transform 
of $D$ on $X_i$, 
for every $i$ {\em{(}}see Theorem $\ref{flip})$. 
Then the sequence terminates after finitely many steps. 
\end{thm}
\begin{proof} 
Suppose that there exists an infinite sequence of 
elementary transformations. 
Let $\Delta$ be the fan corresponding to $X_0$. 
Since the elementary transformations do not change 
one-dimensional cones of $\Delta$, 
there exist numbers 
$k<l$ such that 
the composition $X_k\dashrightarrow X_{k+1}
\dashrightarrow \cdots\dashrightarrow 
X_l$ is an identity. This contradicts the 
negativity in Lemma \ref{neg} below. 
\end{proof}

The following result is easy but very important. 
The proof is well-known. See, for example, \cite[Lemma 3.38]{km}. 

\begin{lem}[the Negativity Lemma]\label{neg}
We consider a commutative diagram 
$$
\begin{matrix} 
&Z&\\
\ \ \ \ \swarrow & & \searrow \ \ \ \ \\
U & \dashrightarrow & \ V \\
{\ \ \ \ \searrow} & \ &  {\swarrow}\ \ \ \ \\
 \ & W &  
\end{matrix}
$$
and $\mathbb Q$-Cartier divisors $D$ and $D'$ 
on $U$ and $V$, respectively, where 
\begin{itemize}
\item[(1)] $f:U\longrightarrow W$ and $g:V\longrightarrow W$ are 
birational morphisms between varieties, 
\item[(2)] $f_*D=g_*D'$, 
\item[(3)] $-D$ is $f$-ample and 
$D'$ is $g$-ample, 
\item[(4)] $\mu:Z\longrightarrow U$, $\nu:Z\longrightarrow V$ are 
common resolutions. 
\end{itemize}
Then $\mu^*D=\nu^*D'+E$, 
where $E$ is an effective $\mathbb Q$-divisor and is exceptional 
over $W$. Moreover, 
if $f$ or $g$ is non-trivial, then $E\ne 0$.  
\end{lem}

\section{On the Zariski decomposition}\label{sec-zari}

In this section, we treat the Zariski decomposition on toric 
varieties. 
Since the MMP works for any divisors, it is obvious that 
the Zariski decomposition holds with no extra assumptions. 

There are many variants of the Zariski decomposition. 
We adopt the following definition 
(cf.~\cite[Definition 7-3-5]{kmm}) here. 

\begin{defn}[the Zariski Decomposition]\label{zari} 
Let $f:X\longrightarrow Y$ be a proper surjective 
morphism of normal varieties. 
An expression $D=P+N$ with $\mathbb R$-Cartier divisors 
$D$, $P$ and 
$N$ on $X$ is called the {\em{Zariski decomposition of $D$ 
relative to $f$ 
in the sense of Cutkosky-Kawamata-Moriwaki}} 
(we write C-K-M for short) if the following conditions are 
satisfied: 
\begin{itemize}
\item[(1)] $P$ is $f$-nef, 
\item[(2)] $N$ is effective, and 
\item[(3)] the natural homomorphisms 
$f_*\mathcal O_X(\llcorner mP\lrcorner)\longrightarrow 
f_*\mathcal O_X(\llcorner mD\lrcorner)$ are bijective for all 
$m\in \mathbb N$. 
\end{itemize}
The divisors $P$ and $N$ are said to be the {\em{positive}} and 
{\em{negative part}} of $D$, respectively.  
\end{defn}
 
\begin{defn}[Pseudo-effective Divisors]\label{pe} 
Let $f:X\longrightarrow Y$ be a projective morphism 
between varieties and $D$ a 
$\mathbb Q$-Cartier divisor on $X$. 
Then $D$ is $f$-{\em{pseudo-effective}} 
if there is an $f$-big (see \cite[Definition 3.22]{km}) 
Cartier divisor $A$ on $X$ such that 
$nD+A$ is $f$-big for every $n\geq 0$. 
\end{defn}

\begin{rem}\label{5333}
Let $f:X\longrightarrow Y$ be a projective morphism 
between varieties and $D$ a 
$\mathbb Q$-Cartier divisor on $X$. 
It is not difficult to see that 
if $D$ is $f$-pseudo-effective then $nD+A$ 
is $f$-big for every $n\geq 0$ 
and any $f$-big Cartier 
divisor $A$ (cf.~\cite[(11.3)]{mori}). 
In particular, if $D$ is an effective divisor on $X$, then 
$D$ is $f$-pseudo-effective. More generally, 
if there exists $m>0$ such that $f_*\mathcal O_X(mD)\ne 0$, 
then $D$ is $f$-pseudo-effective. 
\end{rem}

The following theorem is a slight generalization of 
\cite[Proposition 5]{kawa-zari}. 
Related topics are in \cite[Chapter IV \S 1]{nakayama}. 
Both \cite{kawa-zari} and \cite{nakayama} showed how 
to subdivide a given fan. 

\begin{thm}[{cf.~\cite[Proposition 5]{kawa-zari}}]\label{zarizari}
Let $f:X\longrightarrow Y$ be a projective surjective 
toric 
morphism and $D$ a $\mathbb Q$-Cartier divisor on $X$. 
Assume that $D$ is $f$-pseudo-effective. 
Then there exists a projective birational toric morphism 
$\mu:Z\longrightarrow X$ such that 
$\mu^*D$ has a Zariski decomposition relative 
to $f\circ \mu$ in the sense of {\em{C-K-M}} 
whose positive part is $f\circ \mu$-semi-ample 
$($see \cite[Definition 0-1-4]{kmm}$)$. 
\end{thm}

\begin{rem} 
If $f_*\mathcal O_X(mD)\ne 0$ for some positive integer $m$, 
then 
it is easy to check 
that $(f_k)_{*}\mathcal O_{X_k}(mD_k)\ne 0$ for 
every $k$, where $f_k:X_k\longrightarrow Y$ is as in 
the following proof. 
Therefore, the following proof works without any changes 
even if we replace the assumption that $D$ is 
$f$-pseudo-effective with 
a slightly stronger one that $f_*\mathcal O_X(mD)\ne 0$ 
for some positive integer $m$. 
Thus, we may not have to introduce the 
notion of $f$-pseudo-effective divisors. 
See Corollary \ref{omake} and the proof 
of Corollary \ref{finifini} below.  
\end{rem}

\begin{proof}[Proof of \em{Theorem \ref{zarizari}}] 
By taking a resolution of singularities, 
we may assume that $X$ is non-singular without loss of 
generality. 
Run the MMP on $X$ over $Y$ with respect to 
$D$. 
We obtain a sequence of divisorial contractions and 
elementary transformations over $Y$: 
$$
X=:X_0\dashrightarrow X_1\dashrightarrow X_2\dashrightarrow 
\cdots \dashrightarrow X_k\dashrightarrow X_{k+1}\dashrightarrow 
\cdots. 
$$ 
Since $D_k$ is pseudo-effective over $Y$ for 
every $k$, there exists $l$ such that $D_l$ is nef over 
$Y$ (for the definition of $D_k$, see \ref{mmp}). 
The reader may verify 
that relative pseudo-effectivity of $D$ 
is preserved in each step by Lemma ref{neg} and 
Remark \ref{5333}. 
Take a non-singular quasi-projective toric variety $Z$ with proper 
birational toric morphisms $\mu:Z\longrightarrow X$ and 
$\mu_i:Z\longrightarrow X_i$ for every $0<i\leq l$. 
Then we obtain that 
$\mu^*D=\mu_{l}^* D_l+E$, where $E$ is an effective 
$\mathbb Q$-divisor by the negativity in Lemma \ref{neg}. 
This decomposition is the Zariski 
decomposition of $D$ in the sense of C-K-M. 
\end{proof}

The next corollary is obvious by Theorem \ref{zarizari}. 

\begin{cor}\label{omake} 
Let $f:X\longrightarrow Y$ be a projective surjective toric 
morphism and $D$ a $\mathbb Q$-Cartier divisor on $X$. 
Then $D$ is $f$-pseudo-effective 
if and only if $f_*\mathcal O_X(mD)\ne 0$ 
for some positive integer $m$. 
\end{cor}

\begin{rem}
There exist various generalizations of Theorem \ref{zarizari}. 
We do not pursue such generalizations here. 
For example, the above theorem holds for a 
(not necessarily $\mathbb R$-Cartier) 
$\mathbb R$-divisor $D$, with suitable modifications. 
We leave the details to the reader. 
\end{rem}

The following result is a generalization of Theorem \ref{fini}. 

\begin{cor}[Finitely Generatedness of  Divisorial Algebra II]
\label{finifini}
Let $f:X\longrightarrow Y$ be a proper surjective 
toric morphism and 
$D$ a $($not necessarily $\mathbb Q$-Cartier$)$ 
Weil divisor on $X$. 
Then 
$$
\bigoplus_{m\geq 0}f_*\mathcal O_X(mD)
$$
is a finitely generated $\mathcal O_Y$-algebra. 
\end{cor}
\begin{proof}
By Remarks \ref{qfac} and \ref{qcar}, we may assume that 
$X$ is $\mathbb Q$-factorial. 
Hence, $D$ is $\mathbb Q$-Cartier. 
By replacing $X$ birationally, we may assume that $f$ is 
projective. 
If $f_*\mathcal O_X(mD)=0$ for every $m>0$, 
then the claim is obvious. 
Therefore, we may assume that 
$f_*\mathcal O_X(mD)\ne0$ for some $m>0$, that is, 
$D$ is $f$-pseudo-effective. 
Since by Theorem \ref{zarizari}, there exists a projective 
birational toric morphism $\mu:Z\longrightarrow 
X$ such that $\mu^*D$ 
has a Zariski decomposition with 
$f\circ \mu$-semi-ample positive part, 
$\bigoplus_{m\geq 0}f_*\mathcal O_X(mD)$ is 
finitely generated. 
\end{proof}

\section{Application to hypersurface singularities}\label{appl}

In this section, we apply the toric Mori 
theory to the 
study of singularities. 
We will recover Ishii's results \cite{shihoko}. 
We will work over the complex number field $\mathbb C$, 
throughout this section. 

Let us recall the notion of {\em{non-degenerate}} 
hypersurface singularities quickly. 
For the details, see \cite{shihoko}. 

\begin{defn}[Non-degenerate Polynomials] 
For a polynomial $f=\sum _{\bf m} a_{\bf m}
x^{\bf m}\in \mathbb C[x_0, x_1, \cdots, x_n]$, 
where $x^{\bf m}=x_0^{m_0}x_1^{m_1}\cdots x_n^{m_n}$ for 
${\bf m}=
(m_0, m_1, \cdots, m_n)\in 
{{\mathbb Z}}_{\geq 0}^{n+1}$,  
and a face $\gamma$ of the Newton 
polytope $\Gamma_+(f)$ of $f$, 
denote $\sum _{{\bf m}\in \gamma}a_{\bf m}x^{\bf m}$ 
by $f_\gamma$. 
A polynomial $f$ is said to be 
{\em{non-degenerate}} if for every compact 
face $\gamma$ of $\Gamma _+(f)$ 
the $\partial f_\gamma/ 
\partial x_i \ (i=0, \cdots, n)$ have 
no common zero on $(\mathbb C^\times)^{n+1}$. 
\end{defn}

The following definitions are due to Ishii. See the introduction of \cite{shihoko}. For the definitions of the singularities, 
see Definition \ref{def1}. 

\begin{defn}[Minimal and Canonical Models]\label{model}
Let $(x\in X)$ be a germ of normal singularity on an algebraic 
variety . 
We call a morphism $\varphi:Y\longrightarrow X$ {\em{a minimal 
$($resp.~the canonical$)$ model}} of $(x\in X)$, if 
\begin{itemize}
\item[(1)] $\varphi$ is proper, birational, 
\item[(2)] $Y$ has at most terminal (resp.~
canonical) singularities, and 
\item[(3)] $K_Y$ is $\varphi$-nef 
(resp.~$\varphi$-ample). 
\end{itemize}
\end{defn}

It is obvious that if a canonical model exists, then 
it is unique up to isomorphisms over $X$. 

The next theorem is \cite[Theorem 2.3]{shihoko}. 

\begin{thm}\label{6.3}
Let $X\subset \mathbb C^{n+1}$ be 
a normal hypersurface defined by a non-degenerate polynomial $f$. 
Then $(0\in X)$ has both a minimal model and 
canonical model. 
\end{thm}
\begin{proof}
Take a projective birational toric 
morphism $g:V\longrightarrow \mathbb C^{n+1}$ 
such that $V$ is a non-singular toric variety and 
the proper transform $X'$ of $X$ on $V$ is non-singular 
(see, for example, \cite[Proposition 2.2]{shihoko}). 
Run the MMP over $\mathbb C^{n+1}$ 
with respect to $K_V+X'$. 
Then we obtain $\widetilde{\varphi}:(\widetilde V,\widetilde X)
\longrightarrow \mathbb C^{n+1}$ such that 
$K_{\widetilde V}+\widetilde {X}$ is $\widetilde {\varphi}$-nef. 
We note that the pair $(\widetilde V,\widetilde X)$ 
has canonical singularities and 
$\widetilde V$ has at most terminal singularities. 
So $\widetilde V$ is non-singular in codimension two. 
Thus we obtain 
$K_{\widetilde X}=(K_{\widetilde V}+\widetilde X)|_{\widetilde X}$. 
Therefore, $K_{\widetilde X}$ is nef over $X$. 
It is not difficult to check that $\widetilde X$ has at most 
terminal singularities. 
Hence, this $\widetilde X$ is a minimal model of $(0\in X)$ 
(see Definition 
\ref{model}). 
By using the relative base point free theorem (see, for 
example, \cite[Theorem 
3-1-1 and Remark 3-1-2(1)]{kmm}, or \cite[Theorem 3.24]{km}), 
we obtain the canonical model of $(0\in X)$.  
\end{proof}

\begin{defn}[Dlt and Log-canonical Models]\label{log-model}
Let $(x\in X)$ be as in Definition \ref{model}. 
We call a morphism $\varphi:Y\longrightarrow X$ {\em{a dlt 
$($resp.~the log-canonical$)$ model}} of $(x\in X)$, if 
\begin{itemize}
\item[(1)] $\varphi$ is proper birational, 
\item[(2)] $(Y,E)$ is dlt (resp.~
log-canonical), where $E$ is the reduced exceptional 
divisor of $\varphi$, and 
\item[(3)] $K_Y+E$ is $\varphi$-nef 
(resp.~$\varphi$-ample). 
\end{itemize}
\end{defn}

It is obvious that if a log-canonical model exists, then 
it is unique up to isomorphisms over $X$. 
The notion of dlt models is new.

The next result is a slight generalization of \cite[Theorem 3.1]
{shihoko}. The arguments in the following proof are 
more or less known to the experts of the 
MMP. 

\begin{thm}\label{6.5}
Let $X\subset \mathbb C^{n+1}$ be 
a normal hypersurface defined by a 
non-degenerate polynomial $f$. 
Then $(0\in X)$ has both a minimal model and 
log-canonical model. 
\end{thm}
\begin{proof}
Take a projective birational toric 
morphism $f_0:V_0\longrightarrow \mathbb C^{n+1}$ 
such that $V_0$ is a non-singular toric variety and 
the proper transform $X_0$ of $X$ on $V_0$ is non-singular. 
We may assume that the reduced 
exceptional divisor $E_0$ of $f_0$ 
intersects $X_0$ transversaly, 
that is, $E_0\cup X_0$ is a simple normal crossing divisor 
on $V_0$ (see, for example, \cite[Proposition 2.2]{shihoko}). 
We note that $f_0$ is an isomorphism 
outside $E_0$. 
Run the MMP over $\mathbb C^{n+1}$ 
with respect to $K_{V_0}+X_0+E_0$. 
Then we obtain a sequence of divisorial 
contractions and elementary transformations 
$$
V_0\dashrightarrow V_1
\dashrightarrow V_2\dashrightarrow 
\cdots \dashrightarrow V_k
\dashrightarrow V_{k+1}\dashrightarrow 
\cdots,  
$$ 
and the final object $\widetilde f:\widetilde V\longrightarrow 
\mathbb C^{n+1}$ has the property that 
$K_{\widetilde V}+\widetilde X+\widetilde E$ is $\widetilde f$-nef, 
where $\widetilde X$ is the proper transform of $X_0$ on 
$\widetilde V$ and 
$\widetilde E$ is the reduced $\widetilde f$-exceptional divisor. 
We note that the exceptional locus of 
$V_i\longrightarrow \mathbb 
C^{n+1}$ is of pure codimension one for every $i$ 
since $\mathbb C^{n+1}$ 
is non-singular. 
So, $\widetilde f$ is an isomorphism outside $\widetilde E$. 
Since $(\widetilde V, \widetilde X+\widetilde E)$ is dlt (cf.~
\cite[Corollary 3.44]{km}), 
$\widetilde V$ is non-singular in codimension two 
around $\widetilde X\cap 
\widetilde E$ (cf.~\cite[Corollary 5.55]{km}). 
Thus, we obtain that 
$K_{\widetilde X}+\widetilde E|_{\widetilde X}=
(K_{\widetilde V}+\widetilde X+\widetilde E)|_{\widetilde X}$ 
and $(\widetilde X, \widetilde E|_{\widetilde X})$ is dlt 
(cf.~\cite[Proposition 5.59]{km}). 
We have to check that $\widetilde E|_{\widetilde X}$ 
is a reduced $\widetilde f|_{\widetilde X}$
-exceptional divisor. 
Let $\widetilde E=\sum _i\widetilde E_i$ 
be the irreducible decomposition. 
It is sufficient to show that 
$\widetilde f(\widetilde E_i)\subset \Sing(X)$, 
where $\Sing (X)$ is the singular locus of $X$. 
We write $K_{\widetilde V}+\widetilde X+\sum a_i\widetilde E_i
=\widetilde f^*(K_{\mathbb C^{n+1}}+X)$. 
So, $\sum (1-a_i)\widetilde E_i$ is $\widetilde f$-nef. 
If $\widetilde f(\widetilde E_i)\not\subset \Sing (X)$, 
then $a_i<1$. 
We note that $(\mathbb C^{n+1}, X)$ is plt outside $\Sing (X)$ 
(for the definition of 
{\em{plt}}, see Definition \ref{def1}). 
This implies that $\widetilde f(\widetilde E_i)\subset \Sing (X)$ 
for every $i$ by \cite[Lemma 3.39]{km}. 
Therefore, $(\widetilde X, \widetilde E|_{\widetilde X})$ is 
a dlt model of $(0\in X)$. 
By construction, $\widetilde f|_{\widetilde X}$ is an 
isomorphism outside $\widetilde E|_{\widetilde X}$. 
Since $K_{\widetilde V}+\widetilde X+\widetilde E$ is nef 
over $\mathbb C^{n+1}$, 
we obtain a contraction morphism 
$\widetilde V\longrightarrow V'$ over $\mathbb C^{n+1}$ 
with respect to the divisor 
$K_{\widetilde V}+\widetilde X+\widetilde E$ (cf.~
Theorem \ref{cont}). 
Let $X^{\heartsuit}$ be the normalization of the proper 
transform of $\widetilde X$ on $V'$. 
We put $E^{\heartsuit}:=\mu_*(\widetilde E|_{\widetilde X})$, 
where $\mu:\widetilde X\longrightarrow X^{\heartsuit}$. 
Then it is not difficult to see that 
$K_{X^{\heartsuit}}+E^{\heartsuit}$ 
is ample over $X$, $E^{\heartsuit}$ is the reduced exceptional 
divisor of $X^{\heartsuit}\longrightarrow X$, and 
$K_{\widetilde X}+\widetilde E|_{\widetilde X}
=\mu^*(K_{X^{\heartsuit}}+E^{\heartsuit})$. 
So, $(X^{\heartsuit}, E^{\heartsuit})$ is the required 
log-canonical model of 
$(0\in X)$. 
\end{proof}

\begin{rem}[On Ishii's constructions]\label{shihotan} 
Answering our questions, Ishii informed us that the 
definition of $\overline{E}$ in Claim $3.8$ in the 
proof of \cite[Theorem 3.1]{shihoko} is not correct. 
She told us that $\overline{E}$ should be defined 
as $\nu^*(K_{T_N(\Sigma_2)}+X(\Sigma_2)+E)
-K_{\overline{X}}$, and then all discussions go well 
in the proof of \cite[Theorem 3.1]{shihoko}. 
Though we did not check the proof according to her 
corrected definition of $\overline{E}$, we see 
that our models coincide with her models. 
From now on, we freely use the notation in \cite{shihoko}. 

Let $F$ be the complement of the big torus of $T_N(\Sigma)$. 
Then the pair $(T_N(\Sigma), X(\Sigma))$ (resp.~
$(T_N(\Sigma), X(\Sigma)+F)$) has only canonical 
(resp.~log-canonical) singularities. 
Therefore, by the arguments in 
Claim 2.8 in \cite{shihoko}, 
it is not difficult to see that $m_{\mathbf q}\geq 0$ 
(resp.~$m_{\mathbf q}\geq -1$) for 
every $\mathbf q \in \widetilde \Sigma [1]\setminus \Sigma _0[1]$ 
(resp.~$\mathbf q \in \widetilde \Sigma [1]\setminus \Sigma _2[1]$) 
without the assumption that $D_{\mathbf q}\cap X(\widetilde \Sigma)
\ne \emptyset$ in Claim 2.8 (resp.~Claim 3.5) in \cite{shihoko}. 
Thus, we see that 
$(T_N(\Sigma_0), X(\Sigma_0))$ (resp.~
$(T_N(\Sigma_2), X(\Sigma_2)+E)$) has only canonical (resp.~
log-canonical) singularities. 
Ishii told us that she did not know the notion of singularities of 
pairs when she wrote \cite{shihoko}. 
Therefore, $T_N(\Sigma_0)\simeq 
\Proj_{\mathbb C^{n+1}}
\bigoplus _{m\geq 0}g_*\mathcal O_V(m(K_V+X'))$, 
where $g$, $V$, and $X'$ are as in the proof of Theorem \ref{6.3}, 
and 
$T_N(\Sigma_2)\simeq \Proj
_{\mathbb C^{n+1}}\bigoplus_{m\geq 0}f_{0*}
\mathcal O_X(m(K_{V_0}+X_0+E_0))$, 
where $V_0, X_0, E_0$, and $f_0$ are as 
in the proof of Theorem \ref{6.5}. 
So, it is not difficult to see 
that the models constructed in \cite{shihoko} coincide 
with ours. Details are left to the reader. 
\end{rem}
\ifx\undefined\bysame
\newcommand{\bysame|{leavemode\hbox to3em{\hrulefill}\,}
\fi

\end{document}